\documentclass[a4,twoside,12pt]{article}

\usepackage{a4wide}
\usepackage{amsfonts}
\usepackage{amsmath}
\usepackage{amssymb}
\usepackage{amsthm}
\usepackage[english]{babel}
\usepackage{enumerate}
\usepackage{fontenc}
\usepackage{hyperref}
\usepackage[latin1]{inputenc}
\usepackage{makeidx}
\usepackage[english]{varioref}
\usepackage[all]{xy}
\usepackage{url}

\title{\bf Higgs reductions and numerically flat\\ principal Higgs bundles}
\author{Armando Capasso\\
\\
Universit\`a degli studi di Napoli\,\textquotedblleft Federico II\textquotedblright\\
Scuola Politecnica e delle Scienze di Base\\
Corso Protopisani Nicolangelo 70, Napoli (Italy), C.A.P. 80146\\
\\
E-mail: \texttt{armando.capasso@unina.it}}
\date{}

\DeclareMathOperator{\Ad}{Ad}
\DeclareMathOperator{\Aut}{Aut}
\DeclareMathOperator{\codim}{codim}
\DeclareMathOperator{\End}{End}
\DeclareMathOperator{\Gr}{Gr}
\DeclareMathOperator{\rank}{rank}
\DeclareMathOperator{\tr}{tr}
\DeclareMathOperator{\vol}{Vol}

\theoremstyle{plain}
\newtheorem{theorem}{Theorem}[section]

\newtheorem{corollary}[theorem]{Corollary}
\newtheorem{definition}[theorem]{Definition}
\newtheorem{lemma}[theorem]{Lemma}
\newtheorem{proposition}[theorem]{Proposition}

\theoremstyle{definition}
\newtheorem{example}[theorem]{Example}
\newtheorem{remark}[theorem]{Remark}

\begin{document}

\maketitle

\begin{abstract}
\noindent I consider principal Higgs bundles satisfying a notion of numerical flatness (H-nflatness) that was introduced in~\cite{B:GO:3}. I prove that a principal Higgs bundle $\mathfrak{E}=(E,\varphi)$ is H-nflat is either stable or there exists a Higgs reduction of $\mathfrak{E}$ to a parabolic subgroup $P$ of $G$ such that the principal $L$-bundle $\mathfrak{E}_L$ obtained by extending the reduced Higgs bundle $\mathfrak{E}_P$ to the Levi factor $L$ is H-nflat and stable; and as consequence, $H^{*}(\mathfrak{E},\mathbb{R})$ is isomorphic to the cohomology ring of the associated graded object $\Gr(\mathfrak{E})$ with coefficients in $\mathbb{R}$. Moreover, if $c_2(\Ad(E))$ vanishes then $\mathfrak{E}_L$ is also Hermitian flat and $H^{*}(\Gr(\mathfrak{E}),\mathbb{R})$ is trivial.
\end{abstract}

\tableofcontents

\section{Introduction}
\markboth{Introduction}{Introduction}

\noindent Let $L$ be a line bundle over a complex smooth projective variety $X$ of dimension $n\geq2$: $L$ is said \emph{numerically effective} (\emph{nef}, for short) if nad only if the inequality $L\cdot C\geq0$ holds for any curves $C$ on $X$. Let $V$ be a vector bundle over $X$. $V$ is said \emph{nef} if the line bundle $\mathcal{O}_{\mathbb{P}(V)}(1)$ is nef; where $\mathbb{P}(V)$ is the projectivized bundle. Moreover, $V$ is \emph{numerically flat} (\emph{nflat}, for short) if and only if $V$ and $V^{\vee}$ are nef vector bundles (see~\cite{C:P}).\medskip

\noindent One can also introduce the \emph{Grassmann bundles} associated with $V$: for every integer $s$ such that $0<s<\rank V$, the variety $Gr_s(V)$ is a bundle over $X$, which parameterizes the rank $s$ locally free quotients of $V$. Of course $Gr_1(V)=\mathbb{P}(V)$. Denoting by $\pi_s\colon Gr_s(V)\to X$ the projection, on each variety $Gr_s(V)$ there is a \emph{rank $s$ universal quotient bundle} $Q_s$ of $\pi_s^{*}V$, which turns out to be numerically effective if $V$ is numerically effective.\medskip

\noindent These notions was extend to principal bundle setting in~\cite{B:GO:3} using the bundles $E\left(G_{/\displaystyle P}\right)$; where $G$ is a complex reductive (affine) algebraic group, $P$ is a maximal parabolic subgroup of $G$ and $E$ is a principal $G$-bundle over $X$. Again in the same sense, they parameterize the \emph{reductions of the structure group of $E$ to $P$}; and these are used to provide a notion of \emph{numerically flatness} for principal bundles\footnote{
By~\cite[remark 3.4]{B:GO:3}, when $G=GL(n,\mathbb{C})$ one finds again the Grassmann bundles of the corresponding vector bundle $E$.
} (\emph{nflat}, for short).\medskip

\noindent By~\cite[theorem 6.6]{B:GO:3}, it turns out that the nflat principal $G$-bundles are semistable and can be characterized as follows.
\begin{theorem}\label{th1.1}
Let $E$ be a principal $G$-bundle over $X$. Given the following statements:
\begin{enumerate}[(i)]
\item $E$ is nflat;
\item $E$ is polystable;
\item there exists a reduction $E_P$ of $E$ to a parabolic subgroup $P$ of $G$ such that the principal $L$-bundle $E_P(L)$ is flat and polystable, where $L$ is the Levi factor of $P$;
\end{enumerate}
(i) implies either (ii) or (iii), and (iii) implies (i) and (ii).
\end{theorem}
\noindent Let $\mathfrak{E}=(E,\phi)$ be a principal Higgs $G$-bundle; that is, a principal $G$-bundle equipped with a global section of $\Ad(E)\otimes\Omega^1_X$ such that the composition
\begin{equation*}
[\phi,\phi]:\Ad(E)\xrightarrow{\phi}\Ad(E)\otimes\Omega^1_X\xrightarrow{\phi\otimes Id}\Ad(E)\otimes\Omega^1_X\otimes\Ad(E)\otimes\Omega^1_X\xrightarrow{[\cdot,\cdot]\otimes\pi}\Ad(E)\otimes\Omega^2_X
\end{equation*}
vanishes, where $[\cdot,\cdot]$ is the natural bracket on $\Ad(E)$. In~\cite{B:GO:3} were introduced the \emph{schemes of Higgs reductions} $\mathfrak{R}_P(\mathfrak{E})$ of a principal Higgs $G$-bundle $\mathfrak{E}$ as suitably closed subschemes of $E\left(G_{/\displaystyle P}\right)$, which parameterize the \emph{Higgs reductions of the structure group of $\mathfrak{E}$ to $P$}; and these are used to provide a notion of \emph{numerically flatness} for principal Higgs bundles (\emph{H-nflat}, for short).\medskip

\noindent In this setting, it is known that (i) implies either (ii) or (iii) (cfr.~\cite[remark 6.8]{B:GO:3}). In this paper I prove a partial inversion in three principal steps. Let $\mathfrak{E}$ be a principal Higgs $G$-bundle over $X$:
\begin{enumerate}[a)]
\item if $\mathfrak{E}$ is H-nflat then $\Ad(\mathfrak{E})$ is H-nflat;
\item if $\mathfrak{E}$ is H-nflat. Either it is stable or there exists a Higgs reduction $\mathfrak{E}_P$ of $\mathfrak{E}$ to a parabolic subgroup $P$ of $G$ such that the principal $L$-bundle $\mathfrak{E}_L=\mathfrak{E}_P(L)$ is stable and H-nflat. Moreover, if $c_2(\Ad(\mathfrak{E}))=0$ then $\mathfrak{E}_L$ is also Hermitian flat.
\end{enumerate}

\noindent From now on:
\begin{itemize}
\item $(X,H)$ is an irreducible complex smooth projective variety with a polarization $H$;
\item $\omega$ is a K\"ahler form on $X$ representing the polarization $H$;
\item $G$ is a complex reductive (affine) algebraic group;
\item $\mathfrak{g}$ is the Lie algebra of $G$;
\item $E$ is a principal $G$-bundle over $X$;
\end{itemize}
unless otherwise stated.\medskip

\noindent In the next sections I recall the basic notions about the vector and principal (Higgs) bundles.\medskip

\noindent {\bf Acknowledgement.} This paper is wrote under supervision of professor Ugo Bruzzo (SISSA, Trieste). I am very grateful to him for his useful remarks. This paper was finalized while I was Professor on contract at Scuola Politecnica e delle Scienze di Base of Universit\`a degli studi di Napoli\,\textquotedblleft Federico II\textquotedblright.

\subsection*{Higgs bundles}
\markboth{Higgs bundles}{Higgs bundles}

\begin{definition}\label{def1.1}
A \emph{Higgs sehaf} $\mathfrak{V}$ on $X$ is a coherent sheaf $\mathcal{V}$ on $X$ endowed with a morphism $\phi:\mathcal{V}\to\mathcal{V}\otimes\Omega^1_X$ of $\mathcal{O}_X$-module such that $\phi\wedge\phi=0$. A \emph{Higgs subsheaf} $\mathfrak{W}$ of $\mathcal{V}$ is a \emph{$\phi$-invariant subsheaf $\mathcal{W}$ of $\mathcal{V}$}, that is, $\phi(\mathcal{W})\subseteq\mathcal{W}\otimes\Omega^1_X$. A \emph{Higgs bundle} is a Higgs sheaf $\mathfrak{V}$ such that $\mathcal{V}$ is a locally free $\mathcal{O}_X$-module.
\end{definition}
\noindent Let $\mathfrak{V}=(\mathcal{V},\varphi)$ be a Higgs sheaf over $(X,H)$ with $\rank{\mathcal{V}}>0$. The \emph{slope of $\mathcal{V}$} is defined as
\begin{equation*}
\mu(\mathcal{V})=\frac{c_1(\mathcal{V})\cdot H^{n-1}}{\rank(\mathcal{V})}\equiv\frac{\deg(\mathcal{V})}{\rank(\mathcal{V})}.
\end{equation*}
\begin{definition}\label{def1.2}
A Higgs sheaf $\mathfrak{V}=(V,\phi)$ is \emph{(semi)stable} if $\mathcal{V}$ is torsion free and for any $\varphi$-invariant subsheaf $\mathcal{W}$ with $0<\rank(\mathcal{W})<\rank(\mathcal{V})$ one has
\begin{equation*}
\mu(\mathcal{W})<\mu(V)\,\text{\,(respectively\,}\,\mu(\mathcal{W})\leq\mu(V)\,\text{\,)}.
\end{equation*}
A Higgs sheaf is \emph{polystable} if it is direct sum of stable Higgs sheaves having the same slope.
\end{definition}
\noindent If $\mathfrak{V}=(V,\varphi)$ is a Higgs bundle, and $h$ is an Hermitian metric on $V$, one defines the \emph{Hitchin-Simpson connection} of the pair $(\mathfrak{V},h)$ as
\begin{equation*}
\mathcal{D}_{(h,\varphi)}=D_h+\varphi+\overline{\varphi} 
\end{equation*}
where $D_h$ is the Chern connection of the Hermitian bundle $(V,h)$, and $\overline{\varphi}$ is the metric adjoint of $\varphi$ defined as 
\begin{equation*}
h(s,\varphi(t))=h\left(\overline{\varphi}(s),t\right)
\end{equation*}
for all sections $s,t$ of $V$.
\begin{definition}\label{def1.3}
The \emph{curvature} $\mathcal{R}_{(h,\varphi)}$ and the \emph{mean curvature} $\mathcal{K}_{(h,\varphi)}$ of $\mathcal{D}_{h,\varphi}$ are defined as usual. If $\mathcal{R}_{(h,\varphi)}$ vanishes then $(\mathfrak{V},h)$ is said \emph{Hermitian flat}. If there exists a complex $2$-form on $X$ such that $\mathcal{R}_{(h,\varphi)}=\alpha\cdot I_V$ then $(\mathfrak{V},h)$ is said \emph{projectively flat}.
\end{definition}
\noindent By~\cite[theorem 1.(2)]{S:CT:1992}, the following theorem holds.
\begin{theorem}\label{th1.2}
A Higgs bundle $\mathfrak{V}=(V,\varphi)$ over $X$ is polystable if and only if it admits a Hermitian metric $h$ such that it satisfies the \emph{Hermitian Yang-Mills condition}
\begin{equation*}
\mathcal{K}_{(h,\varphi)}=\frac{2n\pi\mu(V)}{n!\vol(X)}\cdot Id_V\equiv\kappa\cdot Id_V,
\end{equation*}
where $\displaystyle\vol(X)=\int_X\omega^n$.
\end{theorem}
\noindent The following theorem generalize~\cite[theorem IV.4.7 and lemma IV.4.12]{K:S}; the last one will used to conclude the proof of main theorem of this paper.
\begin{theorem}[\bf Bogomolov-L\"ubke's inequality for Higgs bundles]\label{th1.3}
Let $\mathfrak{V}=(V,\varphi)$ be a Higgs bundle over $X$ such that it admits a Hermitian metric $h$ satisfying the Hermitian Yang-Mills condition, then
\begin{equation*}
\int_X\left[2rc_2(V)-(r-1)c_1(V)^2\right]\wedge\omega^{n-2}\geq0
\end{equation*}
and the equality holds if and only if $\displaystyle(\mathfrak{V},h)$ is projectively flat.
\end{theorem}
\noindent {\bf Proof.} On the space of Hermitian endomorphisms $\alpha$ of the Hermitian vector bundle $(V,h)$ I introduce the norm
\begin{equation*}
\|\alpha\|=\max_X\sqrt{\tr(\alpha^2)};
\end{equation*}
by hypothesis $\left\|\mathcal{K}_{(h,\varphi)}\right\|^2=r\kappa^2$ and $\sigma=\tr\left(\mathcal{K}_{(h,\varphi)}\right)=r\kappa$, hence $r\left\|\mathcal{K}_{(h,\varphi)}\right\|^2=\sigma^2$. Let $\{e_i\}$ be a local basis of sections of $V$ and let $\{x^1,\hdots,x^{2n}\}$ be real local coordinates on $X$; one has
\begin{equation*}
\mathcal{R}_{(h,\varphi)}=\frac{1}{2}\sum_{i,j,\alpha,\beta}\left(\mathcal{R}_{(h,\varphi)}\right)^i_{j\alpha\beta}e_i\otimes e_j\otimes dx^i\wedge dx^j,
\end{equation*}
and one puts
\begin{equation*}
\left\|\mathcal{R}_{(h,\varphi)}\right\|^2=\sum_{i,j,\alpha,\beta}\left|\left(\mathcal{R}_{(h,\varphi)}\right)^i_{j\alpha\beta}\right|^2,\,\|\rho\|^2=\sum_{i,\alpha,\beta}\left|\left(\mathcal{R}_{(h,\varphi)}\right)^i_{i\alpha\beta}\right|^2.
\end{equation*}
Using the~\cite[formulas IV.4.2, IV.4.3 and IV.4.6]{K:S}
\begin{gather*}
\begin{split}
\int_X\left[2rc_2(V)-(r-1)c_1(V)^2\right]\wedge\omega^{n-2} &=2r\frac{1}{8\pi^2n(n-1)}\left(\sigma^2-\|\rho\|^2-\left\|\mathcal{K}_{(h,\varphi)}\right\|^2+\left\|\mathcal{R}_{(h,\varphi)}\right\|^2\right)\int_X\omega^n+\\
-(r-1)\frac{1}{4\pi^2n(n-1)}\left(\sigma^2-\|\rho\|^2\right)\int_X\omega^n &=\hdots=\frac{1}{4\pi^2n(n-1)}\left(r\left\|\mathcal{R}_{(h,\varphi)}\right\|^2-\|\rho\|^2\right)\vol(X)\geq0.
\end{split}
\end{gather*}
If the equality holds then by~\cite[formula IV.4.6]{K:S} $\displaystyle\left(\mathcal{R}_{(h,\varphi)}\right)^i_{j\alpha\beta}=\frac{1}{r}\delta^i_j\left(\mathcal{R}_{(h,\varphi)}\right)_{\alpha\beta}$.
\begin{flushright}
(q.e.d.) $\Box$
\end{flushright}
\begin{lemma}\label{lem1.1}
Let $\mathfrak{V}=(V,\varphi)$ be a Higgs bundle over $X$ such that it admits a Hermitian metric $h$ satisfying the Hermitian Yang-Mills condition, $\deg(V)=0$ and $\displaystyle\int_Xc_1(V)^2\wedge\omega^{n-2}=0$. Then
\begin{equation*}
\int_Xc_2(V)\wedge\omega^{n-2}\geq0
\end{equation*}
and the equality holds if and only if the Hitchin-Simpson connection of $(\mathfrak{V},h)$ is Hermitian flat.
\end{lemma}
\noindent {\bf Proof.} By previous theorem, one has the previous inequality. If the equality holds then $\displaystyle\left(\mathcal{R}_{(h,\varphi)}\right)^i_{j\alpha\beta}=\frac{1}{r}\delta^i_j\left(\mathcal{R}_{(h,\varphi)}\right)_{\alpha\beta}$; and by definition, without change the notations by previous proof
\begin{equation*}
\kappa=\frac{2n\pi}{n!\vol(X)}\frac{\deg(V)}{r}=0\Rightarrow\sigma=0
\end{equation*}
that is
\begin{equation*}
0=\int_Xc_2(V)\wedge\omega^{n-2}=\frac{1}{8\pi^2n(n-1)}(r+1)\left\|\mathcal{R}_{(h,\varphi)}\right\|^2\vol(X)\Rightarrow\mathcal{R}_{(h,\varphi)}=0.
\end{equation*}
\begin{flushright}
(q.e.d.) $\Box$
\end{flushright}
\begin{remark}\label{rem1.1}
The theorems~\ref{th1.2},~\ref{th1.3} and the lemma~\ref{lem1.1} hold also for Higgs bundles over compact K\"ahler manifolds.
\begin{flushright}
$\Diamond$
\end{flushright}
\end{remark}

\subsection*{Principal Higgs bundles}
\markboth{Principal Higgs bundles}{Principal Higgs bundles}

\noindent Let $Y$ be a complex variety and let $\rho:G\to\Aut(Y)$ be a representation of $G$ as automorphisms group of $Y$; one defines the vector bundle $E(\rho,Y)=E\times_{\rho}Y$ as the quotient $E\times Y$ under the action of $G$ given by $(e,y)\mapsto\left(eg,\rho\left(g^{-1}\right)y\right)$; in particular, if $Y=\mathfrak{g}$ and $\rho=\Ad$, $E(\Ad,\mathfrak{g})$ is the \emph{adjoint bundle of $E$}, and it is denoted as $\Ad(E)$.\medskip

\noindent Let $\Ad(E)$ be the adjoint bundle of $E$ and let $\Omega_X^1$ be the cotangent sheaf of $X$. If $\varphi$ and $\psi$ are global sections of $\Ad(E)\otimes\Omega_X^1$, one defines a global section $[\varphi,\psi]$ of $\Ad(E)\otimes\Omega_X^2$ by combining the natural bracket $[\cdot,\cdot]:\Ad(E)\otimes\Ad(E)\to\Ad(E)$ and the natural morphism $\pi:\Omega_X^1\otimes\Omega_X^1\to\Omega_X^2$.
\begin{definition}\label{def2.1}
A \emph{Higgs field} $\varphi$ on $E$ is a global section of $\Ad(E)\otimes\Omega_X^1$ such that $[\varphi,\varphi]=0$. A \emph{principal Higgs} $G$\emph{-bundle over} $X$ is a pair $\mathfrak{E}=(E,\varphi)$.
\end{definition}
\noindent Let $A$ be a complex reductive (affine) algebraic group and let $F$ be a principal $A$-bundle over $X$. A \emph{morphism of principal Higgs bundles} $(E,\varphi)\to(F,\psi)$ is a pair $\left(f,f^{\prime}\right)$ such that $f^{\prime}:G\to A$ is an algebraic group homomorphism, $f$ is an $f^{\prime}$-equivariant morphism of bundles over $X$, that is, for any $e\in E$ and $g\in G,\,f(eg)=f(e)f^{\prime}(g)$, and $\left(\widetilde{f}\otimes Id\right)(\varphi)=\psi$; where $\widetilde{f}:\Ad(E)\to\Ad(F)$ is the induced morphism by $f$ between the adjoint bundles given by $\widetilde{f}(e,\alpha)=\left(f(e),f^{\prime}_{*}(\alpha)\right)$, and $f^{\prime}_{*}:\mathfrak{g}\to\mathfrak{a}$ is the morphism of Lie algebras induced by $f^{\prime}$. Let $Y$ be a smooth complex projective variety and let $f:Y\to X$ be a morphism; the \emph{pullback principal Higgs $G$-bundle} $f^{*}\mathfrak{E}$ is the principal $G$-bundle $f^{*}E$ equipped with the Higgs field $f^{*}\varphi$ given by
\begin{equation*}
\varphi\in H^0(X,\Ad(E)\otimes\Omega^1_X)\to H^0\left(Y,\Ad\left(f^{*}E\right)\otimes f^{*}\Omega^1_X\right)\to f^{*}\varphi\in H^0\left(Y,\Ad\left(f^{*}E\right)\otimes\Omega^1_Y\right),
\end{equation*}
where the morphisms are obvious.\medskip

\noindent Let $\mathfrak{a}$ be the Lie algebra of $A$, and let $\lambda:G\to A$ a morphism of algebraic groups. In this case $E(\lambda,A)$ is called \emph{extension (of structure group of)$E$ to $A$}. Where there is not confusion, I shall write $E(\lambda)$ instead of $E(\lambda,A)$. In particular, $G$ acts on $\mathfrak{a}$ via $\Ad_A\circ\lambda$, and one has that the $\mathfrak{a}$-bundle $E\left(\Ad_A\circ\lambda,\mathfrak{a}\right)$ is isomorphic to $\Ad(F)$. Then $\mathfrak{F}=(F,\psi)$ is a principal Higgs $A$-bundle with Higgs field $\psi=\left(\widetilde{f}\times Id\right)(\varphi)$, obtained \emph{extending the structure group $G$ of $\mathfrak{E}$ to $A$}; where $f:E\cong E\times_GG\to F$ is the obvious morphism of principal bundles. In particular, if $\rho:G\to\Aut(\mathbb{V})$ is a linear representation of $G$, the Higgs field of $\mathfrak{E}$ induces a Higgs field on the associated bundle $E\times_{\rho}\mathbb{V}$; moreover, if $\mathbb{V}=\mathfrak{g}$ and $\rho=\Ad$ then $\Ad(\mathfrak{E})$ is the Higgs bundle given by the adjoint bundle $\Ad(E)$ with the induced Higgs field $\Ad(\varphi)$.\medskip

\noindent In a specific way, the construction of \emph{adjoint Higgs bundle} is\,\textquotedblleft invertible\textquotedblright\,as the following proposition states.
\begin{proposition}\label{prop2.1}
Let $E$ be a principal $G$-bundle such that $\Ad(E)$ admits a Higgs field $\phi$. Then there exists a Higgs field $\varphi$ on $E$ (said) \emph{induced by $\phi$}.
\end{proposition}
\noindent {\bf Proof.} Since $\Ad(\Ad(E))\cong\End(\Ad(E))\cong\Ad(E)\otimes\Ad(E)^{\vee}$, by definition $\mathcal{O}_X\xrightarrow{\cdot\phi}\Ad(E)\otimes\Ad(E)^{\vee}\otimes\Omega^1_X$ and $[\phi,\phi]=0\in\Ad(E)\otimes\Ad(E)^{\vee}\otimes\Omega^2_X$. Let $\cdot\varphi:\mathcal{O}_X\xrightarrow{\cdot\phi}\Ad(E)\otimes\Ad(E)^{\vee}\otimes\Omega^1_X\xrightarrow{\pi_1\otimes Id}\Ad(E)\otimes\Omega^1_X$, where $\pi_1\otimes Id$ is the canonical morphism; then
\begin{equation*}
[\varphi,\varphi]=(\pi_1\otimes Id)([\phi,\phi])=0\in\Ad(E)\otimes\Omega^2_X,
\end{equation*}
that is, $\mathfrak{E}=(E,\varphi)$ is a principal Higgs $G$-bundle over $X$.
\begin{flushright}
(q.e.d.) $\Box$
\end{flushright}
\begin{remark}\label{rem2.1}
Considering the principal Higgs $G$-bundle $\mathfrak{E}=(E,\varphi)$, by the isomorphism $\End(E)\cong\Ad(E)\otimes\Ad(E)^{\vee}$, one has the canonical morphism of principal bundles\newline
$f:E\times_{\Ad}\mathfrak{g}\to E\times_{\Ad\otimes\Ad^{\vee}}(\mathfrak{g}\otimes\mathfrak{g}^{\vee})$; by construction $\Ad(\varphi)=\left(\widetilde{f}\times Id\right)(\varphi)$, then the Higgs field induced on $E$ by $\Ad(\varphi)$ is $\varphi$ as well.
\begin{flushright}
$\Diamond$
\end{flushright}
\end{remark}

\section{The scheme of Higgs reductions}
\markboth{The scheme of Higgs reductions}{The scheme of Higgs reductions}

\noindent In this section, I recall the construction of the \emph{scheme of Higgs reductions} of a principal Higgs bundles; and I use it in order to define the notion of \emph{(semi)stable principal Higgs bundle}.\medskip

\noindent Let $K$ be a closed subgroup of $G$ and let $E_{\sigma}$ be a principal $K$-bundle on $X$ together with an injective equivariant principal $K$-bundle morphism $i_{\sigma}:E_{\sigma}\to E$: by definition $E_{\sigma}$ is a \emph{reduction (of the structure group) of $E$ to $K$}. Moreover, these reductions are in one-to-one correspondence with sections $\sigma:X\to E\left(G_{\displaystyle/K}\right)\cong E_{\displaystyle/K}$. By construction there exists an injective morphism of vector bundles $\widetilde{i_{\sigma}}:\Ad(E_{\sigma})\to\Ad(E)$; let $\Pi_{\sigma}:\Ad(E)\otimes\Omega^1_X\to\left(\Ad(E)_{\displaystyle/\Ad(E_{\sigma})}\right)\otimes\Omega^1_X$ the induced projection. Considering the following diagram:
\begin{equation*}
\xymatrix{
0\ar[r] & \Ad(E_{\sigma})\otimes\Omega_X^1\ar[r]^{\widetilde{i_{\sigma}}\otimes Id} & \Ad(E)\otimes\Omega^1_X\ar[r]^(.375){\Pi_{\sigma}} & \left(\Ad(E)_{\displaystyle/\Ad(E_{\sigma})}\right)\otimes\Omega^1_X \ar[r] & 0\\
 & & \mathcal{O}_X\ar[ul]^{\cdot\varphi_{\sigma}}\ar[u]_{\cdot\varphi}
},
\end{equation*}
where $\varphi_{\sigma}\in H^0\left(X,\Ad(E_{\sigma})\otimes\Omega_X^1\right)$ and $[\varphi_{\sigma},\varphi_{\sigma}]=0\in\Ad(E)\otimes\Omega^2_X$.
\begin{definition}[see~\cite{B:GO:3} definition 3.3 and remark 3.4]\label{def2.2}
$\mathfrak{E}_{\sigma}=(E_{\sigma},\varphi_{\sigma})$ is a \emph{Higgs reduction (of structure group of) $\mathfrak{E}$ to $K$} if and only if the previous triangle is commutative, equivalently if and only if $\varphi\in\ker\Pi_{\sigma}$. $\sigma$ is called a \emph{Higgs reduction} of $\mathfrak{E}$.
\end{definition}
\begin{example}\label{ex2.1}
Let $X$ be a compact Riemann surface of genus $g\geq2$, and let $K$ be the canonical bundle of $X$. Considering the \emph{canonical Higgs bundle} $\mathfrak{E}=(E,\varphi)=\left(K^{\frac{1}{2}}\oplus K^{-\frac{1}{2}},\begin{pmatrix}
0 & \omega\\
1 & 0
\end{pmatrix}\right)$ on $X$ of rank $2$, where:
\begin{gather*}
\omega:K^{-\frac{1}{2}}\to K^{\frac{1}{2}}\otimes K\\
1:K^{\frac{1}{2}}\to K^{-\frac{1}{2}}\otimes K\cong K^{\frac{1}{2}}\,\text{\,is the identity section}.
\end{gather*}
By definition~\ref{def2.1} $\varphi$ is a Higgs field; if $\omega\neq0$ then there exists not $\varphi$-invariant line subbundle of $E$. If $\omega=0$ then $K^{-\frac{1}{2}}$ is the unique $\varphi$-invariant line subbundle of $E$, therefore $i:K^{-\frac{1}{2}}\hookrightarrow E$ is an injective equivariant principal $P$-bundle morphism, where $P$ is the parabolic subgroup of $GL(2,\mathbb{C})$; by previous example, $\left(K^{-\frac{1}{2}},0\right)$ is the Higgs reduction of $\mathfrak{E}$ to $P$.
\begin{flushright}
$\triangle$
\end{flushright}
\end{example}
\noindent There is a link between the Higgs reductions of a principal Higgs bundles and the subbundles of the relevant adjoint Higgs bundle, as follows from the next proposition.
\begin{proposition}\label{prop2.2}
Let $\sigma:X\to E_{\displaystyle/K}$ be reduction of $\mathfrak{E}=(E,\varphi)$ to a closed subgroup $K$ of $G$. $\mathfrak{E}_{\sigma}=(E_{\sigma},\varphi_{\sigma})$ is a Higgs reduction of $\mathfrak{E}$ if and only if $\Ad(\mathfrak{E}_{\sigma})$ is a Higgs subbundle of $\Ad(\mathfrak{E})$.
\end{proposition}
\noindent {\bf Proof.} By definition, if $\mathfrak{E}_{\sigma}$ is a Higgs reduction of $\mathfrak{E}$, then $\Ad(\mathfrak{E}_{\sigma})$ is a Higgs subbundle of $\Ad(\mathfrak{E})$. Vice versa, if $\Ad(\mathfrak{E}_{\sigma})$ is a Higgs subbundle of $\Ad(\mathfrak{E})$, by definition $\Ad(\mathfrak{E}_{\sigma})$ is an $\Ad(\varphi_{\sigma})$-invariant subbundle and therefore $\varphi_{\sigma}$ satisfies the definition~\ref{def2.2}; that is, $\mathfrak{E}_{\sigma}$ is a Higgs reduction of $\mathfrak{E}$.
\begin{flushright}
(q.e.d.) $\Box$
\end{flushright}
More in general, the following lemma holds.
\begin{lemma}\label{lem2.1}
Let $\sigma:X\to E_{\displaystyle/K}$ be reduction of $\mathfrak{E}=(E,\varphi)$ to a closed subgroup $K$ of $G$ and let $E_{\sigma}$ the relevant principal $K$-bundle. If $(\Ad(E_{\sigma}),\phi_{\sigma})$ is a Higgs subbundle of $\Ad(\mathfrak{E})=(\Ad(E),\Ad(\varphi))$, then the Higgs field $\varphi_{\sigma}$ induced on $E_{\sigma}$ by $\phi_{\sigma}$ defines a Higgs reduction $\mathfrak{E}_{\sigma}=(E_{\sigma},\varphi_{\sigma})$ $\mathfrak{E}$ to $K$.
\end{lemma}
\noindent {\bf Proof.} By proposition~\ref{prop2.1}, one can consider the principal Higgs $K$-bundle $\mathfrak{E}_{\sigma}=(E_{\sigma},\varphi_{\sigma})$; considering the following diagram
\begin{equation*}
\xymatrix{
\Ad(E_{\sigma})\otimes\Ad(E_{\sigma})\otimes\Omega^1_X\ar[rr]^{\widetilde{j_{\sigma}}}\ar[dd]_{\pi_1\otimes Id} & & \Ad(E)\otimes\Ad(E)\otimes\Omega^1_X\ar[dd]^{\pi_1\otimes Id}\\
& \mathcal{O}_X\ar[ul]^{\cdot\phi_{\sigma}}\ar[ur]^(.4){\cdot\Ad(\varphi)}\ar[dl]^(.4){\cdot\varphi_{\sigma}}\ar[dr]^{\cdot\varphi}\\
\Ad(E_{\sigma})\otimes\Omega^1_X\ar[rr]_{\widetilde{i_{\sigma}}} & & \Ad(E)\otimes\Omega^1_X
}
\end{equation*}
since the square, the upper, right and left triangles are commutative (see also remark~\ref{rem2.1}) then lower triangle is commutative; by definition~\ref{def2.2} $\mathfrak{E}_{\sigma}$ is a Higgs reduction of $\mathfrak{E}$ to $K$.
\begin{flushright}
(q.e.d.) $\Box$
\end{flushright}
Let $\sigma:X\to E_{\displaystyle/K}$ be a Higgs reduction of $\mathfrak{E}$. Given the diagram:
\begin{equation*}
\xymatrix{
\mathcal{O}_{E/K}\ar[d]_{\sigma^{\sharp}}\ar[r]^(.35){\cdot\eta(\varphi)} & T_{E/K,X}\otimes\Omega^1_{E/K}\\
\sigma_{*}\mathcal{O}_X\ar[r]_(.3){\alpha} & \pi_K^{*}\Ad(E)\otimes\Omega^1_{E/K}\ar[u]_{\eta\otimes Id}
}
\end{equation*}
where:
\begin{itemize}
\item $\pi_K$ is the projection of $E_{\displaystyle/K}$ onto $X$;
\item $\gamma$ is the canonical morphism $\Omega^1_X\to\pi^{*}_K\Omega^1_X\to\Omega^1_{E/K}$;
\item $\alpha$ is defined as $(Id\otimes\gamma)\circ\left(\cdot\pi_K^{*}(\varphi)\right)$.
\end{itemize}
Since:
\begin{gather*}
T_{E/K,X}\cong E_{\displaystyle/K}\left(\Ad_K,\mathfrak{g}_{\displaystyle/\mathfrak{k}}\right),\\
\pi_K^{*}\Ad(E)=\pi_K^{*}\left(E(\Ad_G,\mathfrak{g})\right)=E_{\displaystyle/K}(\Ad_K,\mathfrak{g})
\end{gather*}
where $\mathfrak{k}$ is the Lie algebra of $K$, there exists a natural morphism $\eta:\pi_K^{*}\Ad(E)\to T_{E/K,X}$ and $\varphi$ determines a global section $\eta(\varphi)\equiv(\eta\otimes Id)\left(\pi_K^{*}\varphi\right)$ of $T_{E/K,X}\otimes\Omega^1_{E/K}$; from all this: the previous diagram is well defined. Via a straightforward computation on the stalks, one proves that $\eta(\varphi)$ does not depend from $\sigma$, and proves the following statements.
\begin{proposition}[see~\cite{B:GO:3}]\label{prop2.3}
A reduction $\sigma:X\to E_{\displaystyle/K}$ is a Higgs reduction of $\mathfrak{E}=(E,\varphi)$ to $K$ if and only if the scheme-theoretic image of $\sigma$ is contained in the zero locus $V(\eta(\varphi))$ of $\eta(\varphi)$.
\end{proposition}
\begin{corollary}\label{cor2.1}
$V(\eta(\varphi))$ parameterizes the Higgs reductions of $\mathfrak{E}$ to $K$.
\end{corollary}
\begin{definition}[see~\cite{B:GO:3} definition 3.5]\label{def2.3}
$V(\eta(\varphi))$ is called the \emph{scheme of Higgs reduction of $\mathfrak{E}=(E,\varphi)$ to $K$}, and it is denoted as $\mathfrak{R}_K(\mathfrak{E})$.
\end{definition}
\begin{remark}\label{rem2.2}
If $G$ is the general linear group $GL(n,\mathbb{C})$ and $P$ is a maximal parabolic subgroup of $G$, then $G_{\displaystyle/P}$ is the Grassmann variety $Gr_k(\mathbb{C}^n)$ of $k$-dimensional quotient of $\mathbb{C}^n$, for some $k\in\{1,\hdots,n-1\}$, and $\mathfrak{R}_P(\mathfrak{E})$ is the Higgs-Grassmann scheme $\mathfrak{Gr}_k(\mathfrak{V})$ of rank $k$ locally free Higgs quotients of $\mathfrak{V}$, the Higgs bundle corresponding to $\mathfrak{E}$. So a Higgs reduction $\sigma$ of $\mathfrak{E}$ corresponds to a rank $n-k$ Higgs subbundle $\mathfrak{W}$ of $\mathfrak{V}$ (cfr.~\cite{B:GO:1} and~\cite[remark 3.4]{B:GO:3}).
\begin{flushright}
$\Diamond$
\end{flushright}
\end{remark}

\subsection*{(Semi)Stable principal Higgs bundles}
\markboth{(Semi)Stable principal Higgs bundles}{(Semi)Stable principal Higgs bundles}

\noindent Using the scheme of Higgs reduction, I recall the definitions of (semi)stable principal Higgs bundle, given in~\cite{B:GO:3} and some criteria for this.
\begin{definition}[see~\cite{B:GO:3} definition 4.1]\label{def2.5}
A principal Higgs $G$-bundle $\mathfrak{E}$ on $X$ is \emph{(semi)stable} if and only for any parabolic subgroup $P$ of $G$, any open subset $U$ of $X$ such that $\codim(X\setminus U)\geq2$, and any Higgs reduction $\sigma:U\to\mathfrak{R}_P(\mathfrak{E})_{|U}$ of $G$ to $P$ on $U$, one has $\deg\sigma^{*}T_{E/P,X}>0$ (respectively $\deg\sigma^{*}T_{E/P,X}\geq0$).
\end{definition}
\noindent Let $\mathfrak{E}=(E,\varphi)$ be a principal Higgs $G$-bundle over $X$, let $K$ be a closed subgroup of $G$. $E$ is a principal $K$-bundle over $E_{\displaystyle/K}$ via the canonical projection $\pi_K:E\to E_{\displaystyle/K}$. For any character $\chi$ of $K$, one can consider the line bundle $L_{\chi}=E\times_{\chi}\mathbb{C}$ on $E_{\displaystyle/K}$.
\begin{definition}\label{def2.4}
The \emph{slope} of a Higgs reduction $\sigma$ of $\mathfrak{E}$ to $K$ is the homomorphism of groups $\mu_{\sigma}:\chi\in\mathfrak{X}(K)\to\deg\sigma^{*}L_{\chi}\in\mathbb{Q}$; where $\mathfrak{X}(K)$ is the group of characters of $K$.
\end{definition}
\begin{definition}\label{def2.6}
A \emph{dominant character} $\chi$ of $G$ is a character of $G^{\prime}$ (the \emph{semisimple part} of $G$) given by a non negative linear combination of fundamental weights of $\mathfrak{g}^{\prime}$ (the Lie algebra of $G^{\prime}$) such that $\chi$ is trivial on $Z(G)_0$ (the connected component of $1_G$ in $Z(G)$).
\end{definition}
\begin{theorem}[see~\cite{B:GO:3} theorem 4.7]\label{th2.2}
Considered the following conditions:
\begin{enumerate}[a)]
\item for every parabolic subgroup $P$ of $G$ and any dominant character $\chi$ of $P$, the line bundle $L^{\vee}_{\chi}$ restricted to $\mathfrak{R}_P(\mathfrak{E})$ is nef;
\item for every morphism $f:C\to X$, where $C$ is a complex smooth projective curve, the pullback $f^{*}\mathfrak{E}$ is semistable (for short, $\mathfrak{E}$ is Higgs semistable on curves);
\item $\mathfrak{E}$ is semistable and $c_2(\Ad(E))=0\in H^4(X,\mathbb{R})$.
\end{enumerate}
Then the conditions (a) and (b) are equivalent, and they are both implicated by condition (c); and if $\varphi=0$  then the condition (c) is equivalent to conditions (a) and (b).
\end{theorem}

\section{H-nflat (principal) Higgs bundles}
\markboth{H-nflat principal Higgs bundles}{H-nflat principal Higgs bundles}

\noindent In this section, I recall the construction of the \emph{numerically flat principal Higgs bundles} and their main properties.

\subsection*{Universal Higgs quotients}
\markboth{Universal Higgs quotients}{Universal Higgs quotients}

\noindent Let $P$ be a parabolic subgroup of $G$; a connected subgroup $L$ of $G$ is called \emph{Levi factor of $P$} if $P$ is the semi-direct product of $L$ and its \emph{unipotent radical} $R_u(P)$ (see~\cite[definition 11.22 and corollary 14.19]{B:A}). All Levi factors of $P$ are subgroups of $P$ and they are conjugated by elements of $R_u(P)$ (see~\cite[proposition 11.23 and corollary 14.19]{B:A}); they are canonically isomorphic to $P_{\displaystyle/R_u(P)}$ and therefore they are complex connected reductive algebraic groups, whose root systems are, in general, reducible: hence the Levi factor $L$ is decomposable as $L_1\cdot\hdots\cdot L_m$, in according to the decomposition of its root system in irreducible components and any $L_k$ is a simple Lie group (see~\cite[corollary in \S27.5]{H:JE}).\medskip

\noindent Let $\rho:G\to GL(\mathbb{V})$ be a linear representation of $G$ on $\mathbb{V}$, let $\mathbb{W}$ be a vector subspace of $\mathbb{V}$, let $P$ be a maximal parabolic subgroup of $G$ which stabilizes $\mathbb{W}$; then there is an induced action of $P$ on $\mathbb{V}_{\displaystyle/\mathbb{W}}$.
\begin{definition}\label{def3.1}
A factor $L_i$ of the Levi factor $L$ of $P$ is called \emph{standard quotient} of $P$ if $\rho$ maps it injectively in $GL\left(\mathbb{V}_{\displaystyle/\mathbb{W}}\right)$, for some choice of $\rho$ and $\mathbb{W}$.
\end{definition}
\noindent Let $\mathfrak{E}=(E,\varphi)$ be a principal Higgs $G$-bundle over $X$. For any closed subgroup $K$ of $G$, denoted by $E_K$ the principal $K$-bundle $E\to E_{\displaystyle/K}$, one can consider the restriction $i^{*}E_K\equiv i^{*}E\to\mathfrak{R}_K(\mathfrak{E})$ of $E_K$ to $\mathfrak{R}_K(\mathfrak{E})$, where $i$ is the inclusion $\mathfrak{R}_K(\mathfrak{E})\hookrightarrow E_{\displaystyle/K}$. Since $\mathfrak{R}_K(\mathfrak{E})$ parameterizes the Higgs reduction of $\mathfrak{E}$ to $K$, then $\varphi$ induces a Higgs field on the fibres of $i^{*}E_K$ by restriction; that is, $i^{*}E_K$ is a principal Higgs $K$-bundle over $\mathfrak{R}_K(\mathfrak{E})$, and it will be denoted as $\mathfrak{E}_K$. Let $P$ be a (maximal) parabolic subgroup of $G$, let $Q$ be a standard quotient of $P$; let $E_Q$ be the principal $Q$-bundle obtained by extending $E_P$ to $Q$, via the canonical projection $P\to Q$.
\begin{definition}[see~\cite{B:GO:3} definition 5.1]\label{def3.2}
A \emph{universal Higgs quotient $\mathfrak{E}_Q$ of $\mathfrak{E}$} is the restriction of $E_Q$ to the scheme of Higgs reductions $\mathfrak{R}_P(\mathfrak{E})$, equipped with the Higgs field induced by the Higgs fields of $\mathfrak{E}_P$ by extension; here $P$ is maximal parabolic.
\end{definition}
\begin{remark}
\begin{enumerate}[a)]
\item Without change the names from remark~\ref{rem2.2}, a standard quotient $Q$ of $P$ is isomorphic to $GL\left(\mathbb{V}_{\displaystyle/\mathbb{W}}\right)$. In this case, the bundle of local frames of the universal Higgs quotient $\mathfrak{E}_Q$ of $\mathfrak{E}$ is the universal rank $k$ Higgs quotient bundle $\mathfrak{Q}_{k,\mathfrak{V}}$ over $\mathfrak{Gr}_k(\mathfrak{V})$ (cfr.~\cite{B:GO:3} remark 5.4).
\item The construction of the universal Higgs quotients is functorial (see~\cite{B:GO:3} remark 5.3).
\end{enumerate}
\begin{flushright}
$\Diamond$
\end{flushright}
\end{remark}

\subsection*{H-nflat (principal) Higgs bundles}
\markboth{H-nflat (principal) Higgs bundles}{H-nflat Higgs (principal) bundles}

\noindent I recall the definitions of numerically effective (\emph{H-nef}) and numerically flat (\emph{H-nflat}) principal Higgs bundles in the case whose the structure group is an algebraic torus.
\begin{definition}[see~\cite{B:GO:3} definition 5.4]\label{def3.4}
Let $T$ be an algebraic $d$-torus of $G$, where $\dim T=d$, let $\mathfrak{E}=(E,\varphi)$ be a principal Higgs $T$-bundle over $X$.
\begin{enumerate}[a)]
\item $\mathfrak{E}$ is \emph{Higgs-numerically effective} (\emph{H-nef}, for short) if there exists an isomorphism\newline $\lambda:T\to\left(\mathbb{C}^{\times}\right)^d$ such that the bundle $V_{\lambda}$ associated to $E$ via $\lambda$ is nef.
\item $\mathfrak{E}$ is \emph{Higgs-numerically flat} (\emph{H-nflat}, for short) if there exists an isomorphism\newline $\lambda:T\to\left(\mathbb{C}^{\times}\right)^d$ such that the bundle $V_{\lambda}$ associated to $E$ via $\lambda$ is nflat.
\end{enumerate}
\end{definition}
\begin{remark}\label{rem3.1}
\begin{enumerate}[a)]
\item Trivially, the definitions~\ref{def3.4} does not depend from the Higgs field $\varphi$ of $E$.
\item Equivalently, $\mathfrak{E}$ is H-nflat if for \underline{any} isomorphism $\lambda:T\to\left(\mathbb{C}^{\times}\right)^d$ the bundle $V_{\lambda}$ associated to $E$ via $\lambda$ is nflat.
\end{enumerate}
\begin{flushright}
$\Diamond$
\end{flushright}
\end{remark}
\noindent By Structure Theorem for Reductive Groups (see~\cite[proposition in \S14.2]{B:A}) the following sequence
\begin{equation*}
\{1\}\to R(G)\cap G^{\prime}\to R(G)\times G^{\prime}\to G\to\{1\}
\end{equation*}
is short exact; passing to quotients via $G^{\prime}$ (the \emph{derived subgroup} of $G$), one has that $R(G)$ (the \emph{radical} of $G$) is isomorphic to $G_{\displaystyle/G^{\prime}}$; one can define $rad:G\to R(G)$ as the canonical projection of $G$ onto $R(G)$.
\begin{definition}[see~\cite{B:GO:3} definition 5.5]\label{def3.5}
The \emph{radical} of $\mathfrak{E}$ is the principal Higgs $R(G)$-bundle $R(\mathfrak{E})=\mathfrak{E}\times_{rad}R(G)\cong\mathfrak{E}_{\displaystyle/G^{\prime}}$.
\end{definition}
\begin{example}\label{ex3.1}
Let $\mathfrak{V}=(V,\varphi)$ be a rank $n$ Higgs bundle and let $\mathfrak{E}$ be the principal Higgs $GL(n,\mathbb{C})$-bundle of its local frames; then $R(\mathfrak{E})$ is isomorphic to $\mathfrak{E}_{\displaystyle/SL(n,\mathbb{C})}$, the bundle of local frames of the determinant Higgs line bundle $\det(\mathfrak{V})=(\det(V),\det(\varphi))$.
\begin{flushright}
$\triangle$
\end{flushright}
\end{example}
\noindent I recall that the \emph{semi-simple rank $\rank_{ss}(G)$ of $G$} is defined as the rank of the Lie group $G_{\displaystyle/R(G)}$ (see~\cite[definition in \S13.13]{B:A}).
\begin{definition}[see~\cite{B:GO:3} definition 5.7]\label{def3.6}
A principal Higgs $G$-bundle $\mathfrak{E}$ over $X$ is \emph{Higgs-numerically effective} if
\begin{enumerate}[a)]
\item $R(\mathfrak{E})$ is H-nef in according to definition~\ref{def3.4};
\item if $\rank_{ss}(G)>0$, for every maximal parabolic subgroup $P$ of $G$ and every standard quotient $Q$ of $G$, the universal Higgs quotient $\mathfrak{E}_Q$ is H-nef.
\end{enumerate}
Moreover, $\mathfrak{E}$ is \emph{Higgs-numerically flat} if it is H-nef and $R(\mathfrak{E})$ is H-nflat.
\end{definition}
\begin{proposition}\label{prop3.1}
\begin{enumerate}[a)]
\item H-nflat principal Higgs bundles are semistable (\cite[proposition 5.10]{B:GO:3}).
\item If a principal Higgs $G$-bundle $\mathfrak{E}=(E,\varphi)$ is semistable, $c_2(\Ad(E))=0$ and $R(\mathfrak{E})$ is H-nflat, then it is H-nflat (\cite[proposition 5.12]{B:GO:3}).
\end{enumerate}
\end{proposition}
\begin{remark}\label{rem3.2}
\begin{enumerate}[a)]
\item If $G=GL(n,\mathbb{C})$ one has that a Higgs bundle is H-nflat as defined in~\cite{B:GO:1} if and only if the principal Higgs bundle of local frames is H-nflat.
\item A principal Higgs bundle $\mathfrak{E}$ is H-nflat if and only if for any morphism $f:C\to X$, where $C$ is a smooth complex curves, $f^{*}\mathfrak{E}$ is H-nflat (see~\cite[remark 5.8.(ii)]{B:GO:3}).
\end{enumerate}
\begin{flushright}
$\Diamond$
\end{flushright}
\end{remark}

\section{Higgs reductions and\\ H-nflat principal Higgs bundles}
\markboth{Higgs reductions and H-nflat principal Higgs bundles}{Higgs reductions and H-nflat principal Higgs bundles}

\noindent In~\cite{B:GO:3} was proved that a principal Higgs bundle is H-nflat if it satisfies some condition; the inverse holds only if the Higgs field vanishes. Here I prove a partial inverse implication. For this scope, I need recall the following notions on principal Higgs bundles.

\subsection*{Hermitian Yang-Mills-Higgs reductions for\\ H-nflat principal Higgs bundles}
\markboth{Hermitian Yang-Mills-Higgs reductions for H-nflat principal Higgs bundles}{Hermitian Yang-Mills-Higgs reductions for H-nflat principal Higgs bundles}

\noindent Let $K$ be a maximal compact subgroup of $G$; there exists an involution $\iota$ on $\mathfrak{g}$ whose $+1$-eigenspace is the Lie algebra $\mathfrak{k}$ of $K$ (cfr.~\cite[theorem 10.7.2.3]{P:C}), $\iota$ is called \emph{Cartan involution of $G$ relative to $K$}. Let $\mathfrak{E}=(E,\varphi)$ be a principal Higgs $G$-bundle over $X$; one defines the following involution
\begin{equation*}
\forall x\in X,s_x\in\Ad(E)_x,\alpha_x\in\Omega^1_{X,x},\,\iota_x(s_x\otimes\alpha_x)=-\iota_x(s_x)\otimes\overline{\alpha}_x.
\end{equation*}
Let $\sigma$ be a holomorphic reduction of $E$ to $K$: by~\cite[remark 3]{R:S} there exists a unique complex connection $D_{\sigma}$ on $E$ which is compatible with the complex structure on $E$ and with the reduction $\sigma$. From all this, one defines the \emph{Hitchin-Simpson connection $\mathcal{D}_{\sigma,\varphi}$ of $(\mathfrak{E},\sigma)=(E,\varphi,\sigma)$} as the connection $\mathcal{D}_{\sigma,\varphi}=D_{\sigma}+\varphi+\iota(\varphi)$ (see~\cite{B:GO:3}); moreover, by~\cite[corollary in \S3]{R:S}, the previous construction holds even if $K$ is a closed subgroup of $G$. Where there is no confusion about the reduction $\sigma$, I shall write $\mathcal{D}_{\varphi}^K$ instead of $\mathcal{D}_{\sigma,\varphi}$.
\begin{definition}[see~\cite{B:GO:3} definition 6.1]\label{def4.1}
$\mathfrak{E}$ is \emph{Hermitian flat}, if there exists a reduction of its structure group to a maximal compact subgroup $K$ of $G$, such that the corresponding Hitchin-Simpson connection $\mathcal{D}^K_{\varphi}$ is flat.
\end{definition}
\begin{definition}\label{def4.3}
A reduction $\sigma$ of $\mathfrak{E}$ to $K$ is \emph{Hermitian Yang-Mills-Higgs}, if $\sigma$ is a Higgs reduction and there exists an element $\tau$ in the centre of $\mathfrak{g}$ such that the \emph{mean curvature} $\mathcal{K}_{\sigma,\varphi}$ (computed with $\omega$ as usual) is equal to $\tau$.
\end{definition}
\noindent Following~\cite{A:B} and~\cite{B:GO:3}, I posit the next definitions.
\begin{definition}[cfr. definitions~\ref{def2.4},~\ref{def2.6} and~\cite{A:B} definition 3.4]\label{def4.2}
A Higgs reduction $\sigma$ of $\mathfrak{E}$ to a parabolic subgroup $P$ of $G$ is \emph{admissible} if $\mu_{\sigma}(\chi)=0$ for any dominant character $\chi$ of $P$ which vanishes on the centre of $G$.
\end{definition}
\begin{definition}[cfr.~\cite{A:B} definition 3.5]\label{def4.4}
A principal Higgs $G$-bundle $\mathfrak{E}$ is \emph{polystable} if either it is stable or there exist a parabolic subgroup $P$ of $G$ and a Higgs reduction $\sigma$ of $\mathfrak{E}$ to the Levi subgroup $L$ of $P$ such that
\begin{enumerate}[a)]
\item the reduced principal Higgs $L$-bundle $\mathfrak{E}_{\sigma}$ is stable;
\item the principal Higgs $P$-bundle $\mathfrak{E}_{\sigma}(P)$ is an admissible reduction of $\mathfrak{E}$ to $P$.
\end{enumerate}
\end{definition}
\begin{theorem}[cfr.~\cite{B:GO:3} theorem 6.4]\label{th4.2}
A principal Higgs $G$-bundle $\mathfrak{E}$ is \emph{polystable} if and only if it admits an Hermitian Yang-Mills-Higgs reduction to a maximal compact subgroup $K$ of $G$.
\end{theorem}
\begin{lemma}[cfr.~\cite{B:GO:3} proposition 6.5]\label{lem4.1}
A principal Higgs bundle is polystable if and only if its adjoint bundle is polystable.
\end{lemma}
\begin{corollary}\label{cor4.2}
A principal Higgs bundle $\mathfrak{E}$ is stable if its adjoint bundle $\Ad(\mathfrak{E})$ is stable.
\end{corollary}
\noindent {\bf Proof.} By previous lemma, $\mathfrak{E}$ is polystable; if it is not stable, there exists a Higgs reduction $\sigma$ of $\mathfrak{E}$ to the Levi factor $L$ of a parabolic subgroup of $G$, by lemma~\ref{lem2.1} $\Ad(\mathfrak{E}_{\sigma})$ is a Higgs subbundle of $\Ad(\mathfrak{E})$. Since $L$ is a reductive group, $\deg\Ad(\mathfrak{E}_{\sigma})=0$ (see also~\cite[remark 2.2]{R:A}) and this is in contradiction with the stability of $\Ad(\mathfrak{E})$.
\begin{flushright}
(q.e.d.) $\Box$
\end{flushright}
\noindent The following lemma is the cornerstone for the main theorem of this paper.
\begin{lemma}\label{lem4.2}
Let $\mathfrak{E}$ be a H-nflat principal Higgs $G$-bundle over $X$. Then $\Ad(\mathfrak{E})$ is H-nflat.
\end{lemma}
\noindent {\bf Proof.} By remark~\ref{rem3.2}.c, I can assume that $X$ is a curve. By~\cite[lemma 4.7]{A:B} $\Ad(\mathfrak{E})$ is a semistable Higgs bundle; by~\cite[lemma A.7]{B:GO:3} $\Ad(\mathfrak{E})$ is a H-nflat Higgs bundle.
\begin{flushright}
(q.e.d.) $\Box$
\end{flushright}
\noindent I am in position to prove the main theorem of this paper.
\begin{theorem}\label{th4.1}
Let $\mathfrak{E}=(E,\varphi)$ be a H-nflat principal Higgs $G$-bundle. Either it is stable or there exists a Higgs reduction $\sigma$ of $\mathfrak{E}$ to a parabolic subgroup $P$ of $G$, such that the principal $L$-bundle $\widetilde{\mathfrak{E}_{\sigma}}=\mathfrak{E}_{\sigma}(L)$ is H-nflat and stable; where $L$ is the Levi factor of $P$. Moreover, if $c_2(\Ad(E))=0$ then $\widetilde{\mathfrak{E}_{\sigma}}$ is also Hermitian flat.
\end{theorem}
\noindent {\bf Proof.} Without change the notations, by lemma~\ref{lem4.2} $\Ad(\mathfrak{E})$ is H-nflat, so it is semistable. If it is stable, by corollary~\ref{cor4.2} $\mathfrak{E}$ is stable. Otherwise, by~\cite[theorem 3.2]{B:C} it has a Jordan-H\"older filtration
\begin{equation*}
0\subset\mathfrak{V}_0\subset\mathfrak{V}_1\subset\dots\subset\mathfrak{V}_{m-1}\subset\Ad(\mathfrak{E})
\end{equation*}
whose quotients $\mathfrak{Q}_1,\dots,\mathfrak{Q}_{n+1}$ are locally free, stable and H-nflat. The analysis made in the proof of~\cite[proposition 2.10]{A:B} proves that the previous filtration has an odd number of terms, and the middle term $\mathfrak{V}_k$ is isomorphic to a Higgs subbundle $\Ad(\mathfrak{E}_{\sigma})$ of $\Ad(\mathfrak{E})$; by lemma~\ref{lem2.1} $\mathfrak{E}_{\sigma}=(E_{\sigma},\varphi_{\sigma})$ is a Higgs reduction of $\mathfrak{E}$, whose structure group is a parabolic subgroup $P$ of $G$.\medskip

\noindent Let $\widetilde{\mathfrak{E}_{\sigma}}=\left(\widetilde{E_{\sigma}},\widetilde{\varphi_{\sigma}}\right)$ be the principal Higgs $L$-bundle obtained by extending $\mathfrak{E}_{\sigma}$ to $L$. It turns out that $\Ad\left(\widetilde{\mathfrak{E}_{\sigma}}\right)$ is isomorphic to the quotient $\mathfrak{Q}_k$, which is stable and H-nflat. By previous corollary $\widetilde{\mathfrak{E}_{\sigma}}$ is stable.\medskip

\noindent By remark~\ref{rem3.2}.b I can assume that $X$ is a curve. Considering the homomorphism of Lie groups
\begin{equation*}
L\xrightarrow{\alpha}\Aut(\mathfrak{l})\times R(L)\cong\Aut(\mathfrak{l})\times\left(\mathbb{C}^{\times}\right)^d\stackrel{i}{\hookrightarrow}\Aut\left(\mathfrak{l}\oplus\mathbb{C}^d\right)
\end{equation*}
with $\alpha=(\Ad_L,rad)$, gives an injective Lie algebra homomorphism
\begin{equation*}
\alpha_{*}:\mathfrak{l}\to\End(\mathfrak{l})\oplus\mathfrak{r}(L)
\end{equation*}
where $\mathfrak{l}$ and $\mathfrak{r}(L)$ are the Lie algebras of $L$ and $R(L)$, respectively. Let $\mathfrak{V}=\widetilde{\mathfrak{E}_{\sigma}}(\alpha)=\Ad\left(\widetilde{\mathfrak{E}_{\sigma}}\right)\oplus\mathfrak{W}$; because the vector bundle $\mathfrak{W}$ is associated to $R\left(\widetilde{E_{\sigma}}\right)$ via an isomorphism of $R(L)$ to a torus, by~\cite[corollary 2.2]{A:A:B} it is stable and $\deg W=0$; by~\cite[proposition 3.17]{B:GO:2} and definition~\ref{def3.4} $R\left(\widetilde{\mathfrak{E}_{\sigma}}\right)$ is H-nflat, and by proposition~\ref{prop3.1}.b $\widetilde{\mathfrak{E}_{\sigma}}$ is H-nflat.\medskip

\noindent From now on, let $c_2(\Ad(E))=0$; applying~\cite[theorem 3.16]{B:GO:2} and~\cite[theorem 5.2]{B:C}, one has that $\Ad\left(\widetilde{\mathfrak{E}_{\sigma}}\right)$ is Hermitian flat, in particular $c_2\left(\Ad\left(\widetilde{E_{\sigma}}\right)\right)=0$. Since $\widetilde{\mathfrak{E}_{\sigma}}$ is also polystable, by theorem~\ref{th4.2}, it admits a reduction to a maximal compact subgroup of $L$ such that the corresponding Hitchin-Simpson connection $\mathcal{D}$ satisfies the Hermitian Yang-Mills-Higgs condition. Without change the notations, considering $\mathfrak{V}=\Ad\left(\widetilde{\mathfrak{E}_{\sigma}}\right)\oplus\mathfrak{W}$, again one has that $\mathfrak{W}$ is H-nflat, because the underlying vector bundle $W$ is nflat; so, by theorem~\ref{th2.2} $\Delta(W)=0$ and $c_2(W)=0$, these conditions imply that $c_1(W)^2=0$. $\mathcal{D}$ induces Hitchin-Simpson connections on $\Ad\left(\widetilde{\mathfrak{E}_{\sigma}}\right)$ and $\mathfrak{W}$; since both are polystable Higgs bundles, by theorem~\ref{th1.2} and lemma~\ref{lem1.1}, these are Hermitian flat; since $\alpha_{*}$ is injective, $\mathcal{D}$ is Hermitian flat.
\begin{flushright}
(q.e.d.) $\Box$
\end{flushright}
\noindent For simplicity, let $\Gr(\mathfrak{E})$ be the principal Higgs $G$-bundle obtained extending $\widetilde{\mathfrak{E}_{\sigma}}$ to $G$.
\begin{corollary}\label{cor4.1}
If $\mathfrak{E}$ is H-nflat then the cohomology ring $H^{*}(\mathfrak{E},\mathbb{R})$ is isomorphic to $H^{*}(\Gr(\mathfrak{E}),\mathbb{R})$. And if $c_2(\Ad(E))=0$ then this ring is trivial.
\end{corollary}
\noindent {\bf Proof.} In some sense, continuing the previous proof, let $\rho:G\to\Aut(\mathbb{W})$ be a faithful representation of $G$ and let $\mathfrak{V}=\mathfrak{E}(\rho)$ be the Higgs bundle associated to $\mathfrak{E}$ via $\rho$. Miming the proof of~\cite[lemma 3.1]{BI:BU:2}, there exists a flag $0=\mathbb{W}_0\subset\mathbb{W}_1\subset\hdots\subset\mathbb{W}_{m-1}\subset\mathbb{W}_m=\mathbb{W}$ of $P$-modules, such that $R_u(P)$ acts trivially on the quotients. Thus $\rho(P)$ is contained in a parabolic subgroup $Q$ of $\Aut(\mathbb{W})$ and $\rho(L)$ is contained in the Levi subgroup of $Q$. The graded module $\Gr(\mathfrak{V})$ of the filtration of $\mathfrak{V}$ corresponding to $Q$ is isomorphic to the Higgs bundle $(\Gr(\mathfrak{E}))(\rho)$ and, on the other hand, it is homeomorphic to $\mathfrak{V}$. Since this holds for any $\rho$, one has that $\mathfrak{E}$ and $\Gr(\mathfrak{E})$ are homeomorphic; from this follows the claim.
\begin{flushright}
(q.e.d.) $\Box$
\end{flushright}

\end{document}